\documentclass[titlepage,12pt]{article} 
\usepackage{hyperref}
\usepackage{amssymb,amsthm,amsmath} 
\usepackage{mathrsfs}
\usepackage[a4paper]{geometry}
\usepackage{datetime2}
\usepackage[utf8]{inputenc}
\usepackage[italian,english]{babel}

\date{}

\newcommand{\ep}{\varepsilon}
\newcommand{\re}{\mathbb{R}}

\newcommand{\SG}{\mathcal{SG}}
\newcommand{\ft}{\mathscr{F}}

\newtheorem{thm}{Theorem}[section]

\newtheorem{rmk}[thm]{Remark}
\newtheorem{prop}[thm]{Proposition}
\newtheorem{defn}[thm]{Definition}

\title{Global solutions to the Kirchhoff equation with spectral gap data in the energy space}

\author{Marina Ghisi\vspace{1ex}\\ 
{\normalsize Università degli Studi di Pisa} \\
{\normalsize Dipartimento di Matematica}\\ 
{\normalsize PISA (Italy)}\\
{\normalsize e-mail: \texttt{marina.ghisi@unipi.it}}
\and
Massimo Gobbino\vspace{1ex}\\ 
{\normalsize Università degli Studi di Pisa} \\
{\normalsize Dipartimento di Matematica}\\ 
{\normalsize PISA (Italy)}\\  
{\normalsize e-mail: \texttt{massimo.gobbino@unipi.it}}
}

\begin{document}
\maketitle

\begin{abstract}

We prove that the classical hyperbolic Kirchhoff equation admits global-in-time solutions for some classes of initial data in the energy space. We also show that there are enough such solutions so that every initial datum in the energy space is the sum of two initial data for which a global-in-time solution exists.

The proof relies on the notion of spectral gap data, namely initial data whose components vanish for large intervals of frequencies. We do not pass through the linearized equation, because it is not well-posed at this low level of regularity.
\vspace{6ex}

\noindent{\bf Mathematics Subject Classification 2020 (MSC2020):} 
35L90, 35L20, 35L72, 35D30.

		
\vspace{6ex}


\end{abstract}

 
\section{Introduction}

Let $H$ be a real Hilbert space, and let $A$ be a positive self-adjoint operator on $H$ with dense domain $D(A)$. In this paper we consider the evolution equation 
\begin{equation}
u''(t)+m\left(|A^{1/2}u(t)|^{2}\right)Au(t)=0
\label{K:eqn}
\end{equation}
with initial data
\begin{equation}
u(0)=u_{0},
\qquad\qquad
u'(0)=u_{1},
\label{K:data}
\end{equation}
where $m:[0,+\infty)\to[0,+\infty)$ is a suitable nonlinearity. Equation (\ref{K:eqn}) is an abstract version of the hyperbolic partial differential equation
\begin{equation}
u_{tt}(t,x)-m\left(\int_{\Omega}|\nabla u(t,x)|^{2}\,dx\right)\Delta u(t,x)=0,
\nonumber
\end{equation}
with suitable boundary conditions in an open set $\Omega\subseteq\re^{d}$. When $d=1$ or $d=2$, equations of this type are a possible model for small transversal vibrations of elastic strings or membranes. In the case of strings the equation was derived by G.~Kirchhoff in the celebrated monograph~\cite[Section~29.7]{Kirchhoff} after some mathematical simplifications in the full system of local equations of nonlinear elasticity.

Throughout this paper, with the notable exception of Remark~\ref{rmk:m}, we always assume that $m$ satisfies the strict hyperbolicity assumption
\begin{equation}
	m(\sigma)\geq\mu_{1}>0
	\quad\quad
	\forall\sigma\geq 0.
	\label{hp:k-sh}
\end{equation}

When this is the case, it is well-known that problem (\ref{K:eqn})--(\ref{K:data}) admits a unique local solution provided that the nonlinearity $m$ is locally Lipschitz continuous, and initial data $(u_{0},u_{1})$ are in the space $D(A^{3/4})\times D(A^{1/4})$. This result was substantially established by S.~Bernstein in the pioneering paper~\cite{1940-Bernstein}, and then refined by many
authors (see~\cite{1996-TAMS-AroPan} for a modern version).  Existence of global solutions is known in a multitude of different special cases, such as analytic data (\cite{1940-Bernstein,1984-RNM-AroSpa,1992-InvM-DAnSpa,1992-Ferrara-DAnSpa}), quasi-analytic data (\cite{1984-Tokyo-Nishihara,gg:K-Nishihara}), special nonlinearities (\cite{1985-Pohozaev}), dispersive equations and small data (\cite{1980-QAM-GreHu,1993-ARMA-DAnSpa,2005-JDE-Yamazaki,2013-JMPA-MatRuz}), spectral gap data or spectral gap operators (\cite{2005-JDE-Manfrin,2006-JDE-Hirosawa,gg:K-Manfrinosawa,gg:K-Nishihara,2015-NLATMA-Hirosawa}). Nevertheless, none of these results addresses initial data below the regularity threshold of $D(A^{3/4})\times D(A^{1/4})$.

On the other hand, the natural Hamiltonian for (\ref{K:eqn}) is well-defined in the space $D(A^{1/2})\times H$, which for this reason is called the \emph{energy space}. More precisely, if we set
\begin{equation}
M(\sigma):=\int_{0}^{\sigma}m(s)\,ds
\qquad
\forall\sigma\geq 0,
\label{defn:M}
\end{equation}
then a \emph{formal} calculation shows that all solutions to problem (\ref{K:eqn})--(\ref{K:data}) satisfy the energy equality
\begin{equation}
|u'(t)|^{2}+M\left(|A^{1/2}u(t)|^{2}\right)=
|u_{1}|^{2}+M\left(|A^{1/2}u_{0}|^{2}\right)
\label{eqn:energy-eq}
\end{equation}
for every $t$ in the time-interval where they are defined. 

The special role of the space $D(A^{3/4})\times D(A^{1/4})$ becomes clear when we consider the associated linear equation, as we describe in the next paragraph.

\paragraph{\textmd{\textit{Wave equations with time-dependent coefficients}}}

Let us consider the linear wave equation
\begin{equation}
u''(t)+c(t)Au(t)=0,
\label{eqn:linear}
\end{equation}
where $c:[0,+\infty)\to[0,+\infty)$ is a given function (the square of the propagation speed).

This equation was investigated in the seminal paper~\cite{DGCS} and in many subsequent articles. The main result is that problem (\ref{eqn:linear})--(\ref{K:data}) is well-posed in the energy space $D(A^{1/2})\times H$, and more generally in all spaces of the form $D(A^{\alpha+1/2})\times D(A^{\alpha})$, provided that the coefficient $c(t)$ satisfies the strict hyperbolicity assumption
\begin{equation}
c(t)\geq\mu_{1}>0
\qquad
\forall t\geq 0,
\label{hp:c-sh}
\end{equation}
and is (locally) Lipschitz continuous. If the coefficient $c(t)$ is less regular, for example H\"older continuous or just continuous, then much more regularity is required on initial data, for example Gevrey or analytic regularity, otherwise there are situations (and actually many of them, since counterexamples are ``residual'', as shown in~\cite{gg:DGCS-critical}) where solutions do not exist, not even in the sense of distributions.

Now it is immediate to see that equation (\ref{K:eqn}) reduces to (\ref{eqn:linear}) when we set 
\begin{equation}
c(t):=m(|A^{1/2}u(t)|^{2}).
\label{defn:m2c}
\end{equation}

At this point the strict hyperbolicity (\ref{hp:c-sh}) of the coefficient $c(t)$ follows from the strict hyperbolicity (\ref{hp:k-sh}) of the nonlinearity, while the Lipschitz continuity of $c(t)$ is related to the boundedness of its derivative
\begin{eqnarray}
c'(t) & = &
m'\left(|A^{1/2}u(t)|^{2}\right)\langle A^{1/2}u(t),A^{1/2}u'(t)\rangle
\nonumber
\\
& = &
m'\left(|A^{1/2}u(t)|^{2}\right)\langle A^{3/4}u(t),A^{1/4}u'(t)\rangle.
\nonumber
\end{eqnarray}

This is exactly the point where the space $D(A^{3/4})\times D(A^{1/4})$ comes into play. It is the minimal space that guarantees the $C^{1}$ regularity of the function $t\mapsto|A^{1/2}u(t)|^{2}$, even in the case of a wave equation with a constant coefficient. We observe that this regularity is crucial also in order to provide a rigorous justification of the calculation that leads to the energy equality (\ref{eqn:energy-eq}). 

This is also the reason why the space $D(A^{3/4})\times D(A^{1/4})$ has become some sort of insurmountable barrier for Kirchhoff equations. No proof based on a control of $c(t)$ can go any further.

\paragraph{\textmd{\textit{Our contribution}}}

In this paper we surmount the barrier. Indeed, in Theorem~\ref{thm:main} we prove global existence to problem (\ref{K:eqn})--(\ref{K:data}) for some special choices of initial data in the energy space (see Definition~\ref{defn:SG}). On the one hand, these initial data are special because their ``Fourier series'' has infinitely many tails that are very small. On the other hand, these initial data are ``numerous enough'' because any pair $(u_{0},u_{1})$ in the energy space $D(A^{1/2})\times H$ is the sum of two data for which the global solution exists (see Proposition~\ref{prop:sum}). 

We stress that no nontrivial example of local solution was known so far with initial datum in the energy space. We stress also that our solutions exist beyond linearization, and actually despite linearization. What we mean is that in our case the function $t\mapsto|A^{1/2}u(t)|^{2}$ turns out to be just continuous, and not Lipschitz continuous, and therefore if we take one of our solutions $u(t)$ to (\ref{K:eqn}), and we consider the linear equation obtained by freezing the coefficient according to (\ref{defn:m2c}), then there is no guarantee that the linearized equation (\ref{eqn:linear}) admits solutions with non-analytic data (but actually for some mysterious reason it does, at least for the special pair $(u_{0},u_{1})$ we started with). 

Just for completeness, we recall that in~\cite{gg:K-strong-damping} we proved global existence in the energy space for equation (\ref{K:eqn}) with a suitable strong damping. However, the spirit of that result is completely different, because the linear equation (\ref{eqn:linear}) with the same strong damping is well-posed in the energy space also when the coefficient $c(t)$ is less regular. In other words, in that context the Kirchhoff result is a consequence of a linear result, while here it is a purely nonlinear phenomenon.

\paragraph{\textmd{\textit{Overview of the technique}}}

In the proof we exploit the usual ingredients, but we cook them up in an unusual way, so that the final result has a completely different flavor. Actually we exploit only a small subset of the classical ingredients, because many tools are not available in this setting with low regularity.

We can not use higher order energies, and more generally any quantity that contains $A^{\alpha}u'(t)$ for some $\alpha>0$ or $A^{\alpha}u(t)$ for some $\alpha>1/2$. We can not use linearization, which is the key tool in all known fixed point arguments, both on the coefficient, and on the solution. At the end of the day, we can use only the classical Hamiltonian.

What we do is considering the usual finite dimensional approximation \emph{à la Galerkin}, and showing that the resulting sequence $\{u_{n}(t)\}$ is a Cauchy sequence in the energy space. To this end, we write 
\begin{equation}
u_{n}(t)=s_{n,k}(t)+r_{n,k}(t),
\nonumber
\end{equation}
where $s_{n,k}(t)$ and $r_{n,k}(t)$ are the components of $u_{n}(t)$ corresponding to frequencies that are, respectively, smaller and larger than some threshold $k$. Then we set
\begin{equation}
E_{n,k}^{+}(t):=|r_{n,k}'(t)|^{2}+|A^{1/2}r_{n,k}(t)|^{2},
\nonumber
\end{equation}
and we prove that on a fixed time interval it turns out that
\begin{equation}
E_{n,k}^{+}(t)\leq C_{k}\cdot E_{n,k}^{+}(0)
\label{est:Enmk+}
\end{equation}
for a suitable constant $C_{k}$ that does not depend on $n$. This is the key nonlinear point of the proof, for which we introduce an energy that can be thought as a sort of ``modified Hamiltonian restricted to high-frequency components''. We refer to formula (\ref{defn:F}) for the details.

Now we take two functions $u_{n}(t)$ and $u_{m}(t)$, with $n$ and $m$ larger than the threshold $k$, and with a finite dimensional analysis we show that their low-frequency components $s_{n,k}(t)$ and $s_{m,k}(t)$ satisfy
\begin{equation}
|s_{n,k}'(t)-s_{m,k}'(t)|^{2}+|A^{1/2}(s_{n,k}(t)-s_{m,k}(t))|^{2}\leq 
D_{k}\cdot (E_{n,k}^{+}(0)+E_{m,k}^{+}(0))
\label{est:Enmk-}
\end{equation}
for a suitable constant $D_{k}$ that does not depend on $n$ and $m$. From (\ref{est:Enmk+}) and (\ref{est:Enmk-}) we conclude that
\begin{equation}
|u_{n}'(t)-u_{m}'(t)|^{2}+|A^{1/2}(u_{n}(t)-u_{m}(t))|^{2}\leq
(D_{k}+2C_{k})\cdot (E_{n,k}^{+}(0)+E_{m,k}^{+}(0)).
\label{est:un-um}
\end{equation}

The good news is that $E_{n,k}^{+}(0)$ and $E_{m,k}^{+}(0)$ are uniformly small when $k$ is large enough. The bad news is that $C_{k}$ and $D_{k}$ grow rapidly with $k$. This is the point where spectral gap data come into play, because we can define these special data in such a way that the product in the right-hand side of (\ref{est:un-um}) tends to~0 on some sequence of $k$'s.

\paragraph{\textmd{\textit{Conclusions}}}

What does this result tell us? For sure, it dispels the reasonable belief that ``when the linear goes wrong, then also Kirchhoff goes wrong''. Now we know that Kirchhoff equation can admit solutions when its linearization is not well-posed.

The main open problem for Kirchhoff equations (existence of global solution for Sobolev data) remains open. It might happen that problem (\ref{K:eqn})--(\ref{K:data}) has a global solution for all initial data in the energy space, or there might exist solutions, even with initial data in $D(A^{\infty})$, that blow up in a finite time. We do not lean on either side, but we observe that any counterexample would involve a solution that remains bounded in the energy space and blows up in some higher order norm. This is possible only if the energy ``migrates'' from low frequencies to high frequencies. Spectral gap data show that the migration is not possible if there are ``large enough'' holes in the spectrum of initial data. The counter-positive is that energy can migrate only between frequencies that are close enough, and this complicates a lot the search of a counterexample.

\paragraph{\textmd{\textit{Structure of the paper}}}

This paper is organized as follows. In section~\ref{sec:statements} we fix the functional setting, we recall the notion of weak solutions for (\ref{K:eqn}), we introduce spectral gap data, and we state our global existence result. In section~\ref{sec:proofs} prove the main result. Finally, in section~\ref{sec:final} we comment on some natural questions that remain open concerning weak solutions.


\setcounter{equation}{0}
\section{Statements}\label{sec:statements}

Let $H$ be a real Hilbert space, and let $A$ be a self-adjoint operator on $H$ with domain $D(A)$. For the sake of simplicity we assume that there exist an orthonormal basis $\{e_{i}\}_{i\geq 1}$ of $H$, and a nondecreasing sequence $\lambda_{i}\to +\infty$ of positive real numbers such that 
\begin{equation}
Ae_{i}=\lambda_{i}^{2}e_{i}
\qquad
\forall i\geq 1.
\nonumber
\end{equation}

In Remark~\ref{rmk:HA} we discuss a more general setting where the same theory can be developed. Here we just recall that for any such operator the power $A^{\alpha}v$ is defined for every real number $\alpha$, provided that $v$ lies in a suitable subspace $D(A^{\alpha})$.

We are interested in global weak solutions to equation (\ref{K:eqn}) with initial data (\ref{K:data}) in the energy space $D(A^{1/2})\times H$. At this level of regularity we can not expect (\ref{K:eqn}) to be satisfied in the classical sense, not even in the case where $m$ is constant, and therefore we have to rely on the usual notion of weak solution, that we recall below.

\begin{defn}[Global weak solutions]\label{defn:weak-sol}
\begin{em}

A \emph{global weak solution} to problem (\ref{K:eqn})--(\ref{K:data}) is any function
\begin{equation}
u\in C^{0}\left([0,+\infty),D(A^{1/2})\right)\cap C^{1}\left([0,+\infty),H)\right)
\label{hp:reg-u}
\end{equation}
that satisfies (\ref{K:data}) and in addition
\begin{equation}
\int_{0}^{T}\langle u'(t),\varphi'(t)\rangle\,dt=
\int_{0}^{T}m\left(|A^{1/2}u(t)|^{2}\right)\langle A^{1/2}u(t),A^{1/2}\varphi(t)\rangle\,dt
\nonumber
\end{equation}
for every $T>0$ and for every test function $\varphi\in C^{0}\left([0,T],D(A^{1/2})\right)\cap C^{1}\left([0,T],H)\right)$ such that $\varphi(0)=\varphi(T)=0$.

\end{em}
\end{defn}

The notion of weak solution can also be defined in terms of the components of $u$ with respect to the orthonormal basis $\{e_{i}\}$, as follows. Let us consider a function $u$ with the regularity (\ref{hp:reg-u}), let us define the continuous function $c(t)$ according to (\ref{defn:m2c}), and let us consider the components
\begin{equation}
v_{i}(t):=\langle u(t),e_{i}\rangle.
\nonumber
\end{equation}

Then $u$ is a global weak solution to problem (\ref{K:eqn})--(\ref{K:data}) if and only if for every $i\geq 1$ it turns out that $v_{i}\in C^{2}([0,+\infty),\re)$ is the solution to the linear ordinary differential equation
\begin{equation}
v_{i}''(t)+c(t)\lambda_{i}^{2}v_{i}(t)=0
\label{comp:eqn}
\end{equation}
with initial data
\begin{equation}
v_{i}(0)=\langle u_{0},e_{i}\rangle,
\qquad
v_{i}'(0)=\langle u_{1},e_{i}\rangle.
\label{comp:data}
\end{equation}

Now we introduce a special subset of $D(A^{1/2})\times H$, consisting of the so-called \emph{spectral gap data} or \emph{lacunary data}. Analogous spaces are considered in~\cite{2005-JDE-Manfrin,2006-JDE-Hirosawa,gg:K-Manfrinosawa}.

\begin{defn}[Lacunary data]\label{defn:SG}
\begin{em}

For every $(u_{0},u_{1})\in D(A^{1/2})\times H$, and every positive integer $k$, let us set
\begin{equation}
R_{k}(u_{0},u_{1}):=\sum_{i=k+1}^{\infty}\left(\langle u_{1},e_{i}\rangle^{2}+
\lambda_{i}^{2}\langle u_{0},e_{i}\rangle^{2}\right).
\label{defn:Rk}
\end{equation}

We call $\SG(A)$ the set of all pairs $(u_{0},u_{1})\in D(A^{1/2})\times H$ such that
\begin{equation}
R_{k}(u_{0},u_{1})\cdot k^{2}\exp(k\lambda_{k})\leq 1
\qquad
\text{for \emph{infinitely many} integers $k$.}
\label{hp:SG}
\end{equation}

\end{em}
\end{defn}

We stress that (\ref{hp:SG}) is required to be true only for infinitely many values of $k$, and not for every $k$ or eventually, and that the set of integers $k$ for which (\ref{hp:SG}) holds true does depend on $(u_{0},u_{1})$. For this reason $\SG(A)$ is neither a vector space, nor a closed subset of the energy space. Nevertheless, $\SG(A)$ has the following property, which at a first glance looks counterintuitive.

\begin{prop}[Sum property]\label{prop:sum}

For every pair $(u_{0},u_{1})\in D(A^{1/2})\times H$ there exist two pairs $(v_{0},v_{1})$ and $(w_{0},w_{1})$ in $\SG(A)$ such that $u_{0}=v_{0}+w_{0}$ and $u_{1}=v_{1}+w_{1}$.

\end{prop}

The proof of Proposition~\ref{prop:sum} is an exercise on numerical series that we leave to the interested reader (analogous results are proved in~\cite[Theorem~2]{2005-JDE-Manfrin} and~\cite[Proposition~3.2]{gg:K-Manfrinosawa}).

We are now ready to state the main result of this paper.

\begin{thm}[Global solution for spectral gap data in the energy space]\label{thm:main}

Let $H$ be a Hilbert space, and let $A$ be an operator with the properties stated at the beginning of this section. Let $(u_{0},u_{1})\in\SG(A)$ be a pair of lacunary initial data according to Definition~\ref{defn:SG}. Let $m:[0,+\infty)\to[0,+\infty)$ be a function of class $C^{1}$ such that
\begin{equation}
\inf\{m(\sigma):\sigma\geq 0\}>0.
\label{hp:hyperbolicity}
\end{equation}

Then problem (\ref{K:eqn})--(\ref{K:data}) admits at least one global weak solution in the sense of Definition~\ref{defn:weak-sol}, and this solution satisfies the energy equality (\ref{eqn:energy-eq}) for every $t\geq 0$.

\end{thm}

In the previous statement we were quite generous with the assumptions, because in that framework we can present a clear proof that conveys the key ideas without unnecessary technicalities. Now we discuss some possible extensions.

\begin{rmk}[Assumptions on the nonlinearity]\label{rmk:m}
\begin{em}

The assumptions on the nonlinearity $m$ can be relaxed a little with standard changes in the proof. More precisely, the $C^{1}$ regularity can be weakened to local Lipschitz continuity, and the strict hyperbolicity assumption (\ref{hp:hyperbolicity}) can be weakened by just asking that $m(\sigma)>0$ for every $\sigma\geq 0$ and
\begin{equation}
\int_{0}^{+\infty}m(\sigma)\,d\sigma=+\infty.
\nonumber
\end{equation}

\end{em}
\end{rmk}

\begin{rmk}[Variants of lacunary data]\label{rmk:SG}
\begin{em}

In the definition of $\SG(A)$ we can replace the weight $k^{2}\exp(k\lambda_{k})$ by $\alpha_{k}\exp(\beta_{k}\lambda_{k})$, where $\{\alpha_{k}\}$ and $\{\beta_{k}\}$ are two sequences of real numbers such that $\alpha_{k}\to +\infty$ and $\beta_{k}\to +\infty$ (these two sequences might also depend on $(u_{0},u_{1})$). We can also replace the~1 in the right-hand side of (\ref{hp:SG}) with a constant $C$, that is also allowed to depend on $(u_{0},u_{1})$.  Note that in this way the set $\SG(A)$ becomes invariant by homothety.

We also observe that, if the eigenvalues of $A$ grow fast enough, for example $\lambda_{k+1}\geq\exp(k^{2}\lambda_{k})$ for infinitely many $k$'s, then $D(A^{\alpha+1/2})\times D(A^{\alpha})\subseteq\SG(A)$ for every $\alpha>0$.

\end{em}
\end{rmk}

\begin{rmk}[Assumptions on the operator]\label{rmk:HA}
\begin{em}

Concerning $H$ and $A$, we do not need the operator to be positive, but just nonnegative, and we do not need its spectrum to be discrete. We just need that $H$ is a real Hilbert space, and $A$ is unitary equivalent to a multiplication operator on some $L^{2}$ space. More precisely, this means that there exist a measure space $(\mathcal{M},\mu)$, a linear bijective isometry $\ft:H\to L^{2}(\mathcal{M},\mu)$, and a nonnegative measurable function $\lambda:\mathcal{M}\to[0,+\infty)$ such that
\begin{equation}
v\in D(A)
\quad\Longleftrightarrow\quad
\lambda(\xi)^{2}\cdot[\ft(v)](\xi)\in L^{2}(\mathcal{M},\mu),
\nonumber
\end{equation}
and such that for every $v\in D(A)$ it turns out that
\begin{equation}
[\ft(Av)](\xi)=\lambda(\xi)^{2}\cdot[\ft(v)](\xi)
\quad
\text{for $\mu$-almost every }\xi\in \mathcal{M}.
\nonumber
\end{equation}

In other words, we can identify $H$ with $L^{2}(\mathcal{\mathcal{M}},\mu)$ through $\ft$, which plays the role of a generalized Fourier transform, and under this identification the operator $A$ becomes the multiplication operator by $\lambda(\xi)^{2}$. 

In this context we say that $(u_{0},u_{1})\in\SG(A)$ if and only if there exists a sequence of real numbers $\ell_{n}\to +\infty$ such that
\begin{equation}
n^{2}\exp(n\ell_{n})\cdot
\int_{\{\xi\in \mathcal{M}\,:\,\lambda(\xi)>\ell_{n}\}}
\left(\widehat{u}_{1}(\xi)^{2}+\lambda(\xi)^{2}\widehat{u}_{0}(\xi)^{2}\right)\,d\mu(\xi)\leq 1,
\nonumber
\end{equation}
where as usual we set $\widehat{u}_{0}:=\ft(u_{0})$ and $\widehat{u}_{1}:=\ft(u_{1})$.

\end{em}
\end{rmk}

\bigskip


\setcounter{equation}{0}
\section{Proof of Theorem~\ref{thm:main}}\label{sec:proofs}

We follow a usual Galerkin approximation scheme. To this end, for every positive integer $n$ we consider the solution $u_{n}(t)$ to equation (\ref{K:eqn})
with initial data
\begin{equation}
u_{n}(0)=\sum_{i=1}^{n}\langle u_{0},e_{i}\rangle e_{i},
\qquad\qquad
u_{n}'(0)=\sum_{i=1}^{n}\langle u_{1},e_{i}\rangle e_{i}.
\label{Galerkin:data}
\end{equation}

We observe that the initial data of $u_{n}(t)$ are the projections of the initial data (\ref{K:data}) onto the $n$-dimensional subspace of $H$ generated by $e_{1}$, \ldots, $e_{n}$. Since this finite dimensional subspace of $H$ is $A$-invariant, equation (\ref{K:eqn}) reduces to a system of finitely many ordinary differential equations, and therefore problem (\ref{K:eqn})--(\ref{Galerkin:data}) admits exactly one solution $u_{n}(t)$. This solution, defined for every positive time, is actually a strong solution, and even better it satisfies $u_{n}\in C^{3}([0,+\infty),D(A^{\alpha}))$ for every real number $\alpha$ (the only obstruction to further time regularity of $u_{n}(t)$ is the regularity of $m(\sigma)$). As a consequence, $u_{n}(t)$ satisfies also the energy equality
\begin{equation}
|u_{n}'(t)|^{2}+M\left(|A^{1/2}u_{n}(t)|^{2}\right)=
|u_{n}'(0)|^{2}+M\left(|A^{1/2}u_{n}(0)|^{2}\right)
\qquad
\forall t\geq 0.
\label{eqn:EE-un}
\end{equation}

The key claim is that $\{u_{n}(t)\}$ is a Cauchy sequence in the space 
$$C^{0}([0,T],D(A^{1/2}))\cap C^{1}([0,T],H)$$ 
for every $T>0$. 

If this is true, then there exists $u\in C^{0}([0,+\infty),D(A^{1/2}))\cap C^{1}([0,+\infty),H)$ such that
\begin{equation}
\lim_{n\to +\infty}\sup_{t\in[0,T]}\left(|u'(t)-u_{n}'(t)|^{2}+|A^{1/2}(u(t)-u_{n}(t))|^{2}\right)=0
\qquad
\forall T>0.
\nonumber
\end{equation}

This uniform convergence on bounded intervals is enough to pass to the limit both in the notion of weak solution, from which we deduce that $u(t)$ is a global weak solution to problem (\ref{K:eqn})--(\ref{K:data}) in the sense of Definition~\ref{defn:weak-sol}, and in the energy equality (\ref{eqn:EE-un}), from which we obtain (\ref{eqn:energy-eq}).

In order to prove that $\{u_{n}(t)\}$ is a Cauchy sequence it is enough to show that, for every $T>0$ and every $\ep>0$, there exists a positive integer $k$ such that
\begin{equation}
\sup_{t\in[0,T]}\left(|u_{m}'(t)-u_{n}'(t)|^{2}+|A^{1/2}(u_{m}(t)-u_{n}(t))|^{2}\right)\leq\ep
\qquad
\forall n>k,
\quad
\forall m>k.
\label{th:Cauchy}
\end{equation}

The rest of the proof is devoted to the proof of this claim. To this end, we introduce some notation.

\paragraph{\textmd{\textit{Notation and basic uniform estimates}}}

Let us set 
\begin{equation}
H_{0}:=|u_{1}|^{2}+M(|A^{1/2}u_{0}|^{2}).
\label{defn:E0}
\end{equation}

We observe that the right-hand side of (\ref{eqn:EE-un}) is less than or equal to $H_{0}$, and therefore \begin{equation}
|u_{n}'(t)|^{2}+M\left(|A^{1/2}u_{n}(t)|^{2}\right)\leq
H_{0}
\qquad
\forall t\geq 0,
\quad
\forall n\geq 1.
\label{est:un-H0}
\end{equation}

If $\mu_{1}$ denotes the infimum in (\ref{hp:hyperbolicity}), then from (\ref{defn:M}) we obtain that $M(\sigma)\geq\mu_{1}\sigma$ for every $\sigma\geq 0$, and therefore from (\ref{est:un-H0}) we deduce that there exists a constant $H_{1}$ such that
\begin{equation}
|A^{1/2}u_{n}(t)|^{2}\leq H_{1}
\qquad
\forall t\geq 0,
\quad
\forall n\geq 1.
\nonumber
\end{equation}

Due to this bound, in the sequel we can assume, without loss of generality, that $m(\sigma)$ is also bounded from above and Lipschitz continuous, and more precisely that
\begin{equation}
0<\mu_{1}\leq m(\sigma)\leq\mu_{2}
\qquad
\forall\sigma\geq 0
\label{bound:m}
\end{equation}
and
\begin{equation}
|m(\sigma_{2})-m(\sigma_{1})|\leq L|\sigma_{2}-\sigma_{1}|
\qquad
\forall(\sigma_{1},\sigma_{2})\in[0,+\infty)^{2}
\label{hp:m-lip}
\end{equation}
for suitable positive real numbers $\mu_{1}$, $\mu_{2}$ and $L$. 

For every pair of positive integers $n>k$, we write $u_{n}(t)$ as
\begin{equation}
u_{n}(t)=s_{n,k}(t)+r_{n,k}(t),
\nonumber
\end{equation}
where
\begin{equation}
s_{n,k}(t):=\sum_{i=1}^{k}\langle u_{n}(t),e_{i}\rangle e_{i}
\nonumber
\end{equation}
is the component of $u_{n}(t)$ with respect to the subspace of $H$ generated by $e_{1}$, \ldots, $e_{k}$ (the small frequencies), and
\begin{equation}
r_{n,k}(t):=\sum_{i=k+1}^{n}\langle u_{n}(t),e_{i}\rangle e_{i}
\nonumber
\end{equation}
is the remainder, corresponding to the component of $u_{n}(t)$ with respect to the subspace of $H$ generated by $e_{k+1}$, \ldots, $e_{n}$ (the high frequencies). 

It is easy to check that $s_{n,k}(t)$ solves equation
\begin{equation}
s_{n,k}''(t)+m\left(|A^{1/2}s_{n,k}(t)|^{2}+|A^{1/2}r_{n,k}(t)|^{2}\right)As_{n,k}(t)=0
\nonumber
\end{equation}
with initial data
\begin{equation}
s_{n,k}(0)=\sum_{i=1}^{k}\langle u_{0},e_{i}\rangle e_{i},
\qquad
s_{n,k}'(0)=\sum_{i=1}^{k}\langle u_{1},e_{i}\rangle e_{i},
\nonumber
\end{equation}
while $r_{n,k}(t)$ solves equation
\begin{equation}
r_{n,k}''(t)+m\left(|A^{1/2}s_{n,k}(t)|^{2}+|A^{1/2}r_{n,k}(t)|^{2}\right)Ar_{n,k}(t)=0
\label{eqn:rnk}
\end{equation}
with initial data
\begin{equation}
r_{n,k}(0)=\sum_{i=k+1}^{n}\langle u_{0},e_{i}\rangle e_{i},
\qquad
r_{n,k}'(0)=\sum_{i=k+1}^{n}\langle u_{1},e_{i}\rangle e_{i},
\nonumber
\end{equation}

For the sake of shortness we introduce the function
\begin{equation}
\varphi_{n,k}(t):=|A^{1/2}s_{n,k}(t)|^{2}.
\label{defn:phink}
\end{equation}
and the constants
\begin{equation}
\nu_{1}:=\min\{1,\mu_{1}\},
\qquad\qquad
\nu_{2}:=\max\{1,\mu_{2}\}.
\nonumber
\end{equation}

\paragraph{\textmd{\textit{Uniform bound on low-frequency components}}}

We show that for every pair of positive integers $n>k$ it turns out that
\begin{equation}
|s_{n,k}'(t)|^{2}+|A^{1/2}s_{n,k}(t)|^{2}\leq
\frac{H_{0}}{\nu_{1}}
\qquad
\forall t\geq 0
\label{th:snk}
\end{equation}
and
\begin{equation}
|\varphi_{n,k}'(t)|\leq
\frac{H_{0}}{\nu_{1}}\cdot\lambda_{k}
\qquad
\forall t\geq 0.
\label{th:phink}
\end{equation}

To this end, since $s_{n,k}(t)$ is a projection of $u_{n}(t)$ onto a subspace, from the energy inequality (\ref{est:un-H0}) we deduce that
\begin{equation}
|s_{n,k}'(t)|^{2}+M\left(|A^{1/2}s_{n,k}(t)|^{2}\right)\leq
|u_{n}'(t)|^{2}+M\left(|A^{1/2}u_{n}(t)|^{2}\right)\leq
H_{0}.
\nonumber
\end{equation}

Recalling that $M(\sigma)\geq\mu_{1}\sigma$ for every $\sigma\geq 0$, this is enough to deduce (\ref{th:snk}). Now we observe that
\begin{equation}
\varphi_{n,k}'(t)=2\langle As_{n,k}(t),s_{n,k}'(t)\rangle
\qquad
\forall t\geq 0.
\nonumber
\end{equation}

Since $s_{n,k}(t)$ belongs to the subspace of $H$ generated by $e_{1}$, \ldots, $e_{k}$, and $\lambda_{k}$ is the largest eigenvalue of $A$ in this subspace, we conclude that
\begin{eqnarray*}
|\varphi_{n,k}'(t)| & \leq &
2|As_{n,k}(t)|\cdot|s_{n,k}'(t)|
\\
& \leq &
\lambda_{k}\cdot 2|A^{1/2}s_{n,k}(t)|\cdot|s_{n,k}'(t)|
\\
& \leq &
\lambda_{k}\left(|s_{n,k}'(t)|^{2}+|A^{1/2}s_{n,k}(t)|^{2}\right).
\end{eqnarray*}

At this point (\ref{th:phink}) follows from (\ref{th:snk}).

\paragraph{\textmd{\textit{Uniform smallness of high-frequency components}}}

We show that for every pair of positive integers $n>k$ it turns out that
\begin{eqnarray}
|r_{n,k}'(t)|^{2}+|A^{1/2}r_{n,k}(t)|^{2} & \leq &
R_{k}(u_{0},u_{1})\cdot L_{k}^{+}(t)
\qquad
\forall t\geq 0,
\label{th:rnk}
\end{eqnarray}
where $R_{k}(u_{0},u_{1})$ is defined by (\ref{defn:Rk}) and
\begin{equation}
L_{k}^{+}(t):=\frac{\nu_{2}}{\nu_{1}}\cdot
\exp\left(\frac{LH_{0}\lambda_{k}t}{\nu_{1}^{2}}\right).
\label{defn:L+}
\end{equation}

To this end, we consider the function $\varphi_{n,k}(t)$ defined by (\ref{defn:phink}), and we introduce the energy
\begin{equation}
F(t):=|r_{n,k}'(t)|^{2}+
M\left(\varphi_{n,k}(t)+|A^{1/2}r_{n,k}(t)|^{2}\right)-M\left(\varphi_{n,k}(t)\right).
\label{defn:F}
\end{equation}

From the bounds in (\ref{bound:m}) we deduce that
\begin{equation}
\mu_{1}|A^{1/2}r_{n,k}(t)|^{2}\leq
M\left(\varphi_{n,k}(t)+|A^{1/2}r_{n,k}(t)|^{2}\right)-M\left(\varphi_{n,k}(t)\right)\leq
\mu_{2}|A^{1/2}r_{n,k}(t)|^{2},	
\nonumber
\end{equation}
and therefore
\begin{equation}
\nu_{1}\left(|r_{n,k}'(t)|^{2}+|A^{1/2}r_{n,k}(t)|^{2}\right)
\leq F(t)
\leq
\nu_{2}\left(|r_{n,k}'(t)|^{2}+|A^{1/2}r_{n,k}(t)|^{2}\right).
\label{th:F-equiv}
\end{equation}

Due to the time regularity of $u_{n}(t)$, and hence also of $r_{n,k}(t)$, we can compute the time derivative of $F(t)$. Exploiting (\ref{eqn:rnk}) we obtain that
\begin{equation}
F'(t)=
\varphi_{n,k}'(t)
\left\{m\left(\varphi_{n,k}(t)+|A^{1/2}r_{n,k}(t)|^{2}\right)-m\left(\varphi_{n,k}(t)\right)\right\}.
\nonumber
\end{equation}

Thus from (\ref{th:phink}) and the Lipschitz continuity of $m$ we deduce that
\begin{equation}
F'(t)\leq
\frac{H_{0}\lambda_{k}}{\nu_{1}}\cdot L|A^{1/2}r_{n,k}(t)|^{2}\leq
\frac{LH_{0}\lambda_{k}}{\nu_{1}^{2}}F(t).
\nonumber
\end{equation}

Integrating this differential inequality we obtain that
\begin{equation}
F(t)\leq F(0)
\exp\left(\frac{LH_{0}\lambda_{k}t}{\nu_{1}^{2}}\right)
\qquad
\forall t\geq 0.
\label{est:F}
\end{equation}

Finally, we recall that
\begin{eqnarray}
F(0) & \leq & 
\nu_{2}\left(|r_{n,k}'(0)|^{2}+|A^{1/2}r_{n,k}(0)|^{2}\right)
\nonumber
\\
& = &
\nu_{2}\sum_{i=k+1}^{n}\left(
\langle u_{1},e_{i}\rangle^{2}+\lambda_{i}^{2}\langle u_{0},e_{i}\rangle^{2}
\right)
\nonumber
\\
& \leq &
\nu_{2}R_{k}(u_{0},u_{1})
\label{est:F(0)}
\end{eqnarray}

At this point (\ref{th:rnk}) follows from (\ref{est:F}), (\ref{est:F(0)}), and the estimate from below in (\ref{th:F-equiv}).

\paragraph{\textmd{\textit{Smallness of differences between low-frequency components}}}

We show that for every triple of positive integers $(n,m,k)$, with $n>k$ and $m>k$, it turns out that
\begin{equation}
|s_{n,k}'(t)-s_{m,k}'(t)|^{2}+\left|A^{1/2}(s_{n,k}(t)-s_{m,k}(t))\right|^{2}\leq 
R_{k}(u_{0},u_{1})^{2}\cdot L_{k}^{-}(t)
\qquad
\forall t\geq 0,
\label{th:s-s}
\end{equation}
where
\begin{equation}
L_{k}^{-}(t):=2L\left(\frac{\nu_{2}}{\nu_{1}}\right)^{2}
\exp\left(\left(1+\mu_{2}+\frac{4LH_{0}}{\nu_{1}^{2}}\right)\lambda_{k}t\right).
\label{defn:L-}
\end{equation}

To this end, we introduce the difference
\begin{equation}
\rho_{n,m,k}(t):=s_{n,k}(t)-s_{m,k}(t),
\nonumber
\end{equation}
and we observe that it is a solution to equation
\begin{equation}
\rho_{n,m,k}''(t)+A\rho_{n,m,k}(t)=\psi_{n,k}(t)-\psi_{m,k}(t)
\label{eqn:rhonmk}
\end{equation}
with initial data
\begin{equation}
\rho_{n,m,k}(0)=0,
\qquad
\rho_{n,m,k}'(0)=0,
\nonumber
\end{equation}
where
\begin{gather*}
\psi_{n,k}(t):=
\left\{1-m\left(|A^{1/2}u_{n}(t)|^{2}\right)\right\}As_{n,k}(t),
\\
\psi_{m,k}(t):=
\left\{1-m\left(|A^{1/2}u_{m}(t)|^{2}\right)\right\}As_{m,k}(t).
\end{gather*}

Now we introduce the energy
\begin{equation}
G(t):=|\rho_{n,m,k}'(t)|^{2}+|A^{1/2}\rho_{n,m,k}(t)|^{2}.
\nonumber
\end{equation}

Due to the time regularity of $u_{n}(t)$ and $u_{m}(t)$, and hence also of $\rho_{n,m,k}(t)$, we can compute the time derivative of $G(t)$. Exploiting (\ref{eqn:rhonmk}) we obtain that
\begin{equation}
G'(t)=
2\langle\psi_{n,k}(t)-\psi_{m,k}(t),\rho_{n,m,k}'(t)\rangle\leq
2|\psi_{n,k}(t)-\psi_{m,k}(t)|\cdot|\rho_{n,m,k}'(t)|.
\nonumber
\end{equation}

In order to estimate the first factor, we set for simplicity
\begin{equation}
L_{1}:=\frac{H_{0}}{\nu_{1}},
\nonumber
\end{equation}
and we observe that
\begin{eqnarray}
\psi_{n,k}(t)-\psi_{m,k}(t) & = &
\left\{1-m\left(|A^{1/2}u_{n}(t)|^{2}\right)\right\}A\rho_{n,m,k}(t)
\nonumber
\\[1ex]
& &
+\left\{m\left(|A^{1/2}u_{m}(t)|^{2}\right)-m\left(|A^{1/2}u_{n}(t)|^{2}\right)\right\}As_{m,k}(t).
\label{psi-split}
\end{eqnarray}

Now we estimate the two terms separately. For the first term, we exploit the bound from above in (\ref{bound:m}), and the fact that $\rho_{n,m,k}(t)$ lies in the subspace of $H$ generated by $e_{1}$, \ldots, $e_{k}$, where $\lambda_{k}$ is the largest eigenvalue of $A$. We obtain that
\begin{equation}
\left|1-m\left(|A^{1/2}u_{n}(t)|^{2}\right)\right|\cdot|A\rho_{n,m,k}(t)|\leq
(1+\mu_{2})\lambda_{k}|A^{1/2}\rho_{n,m,k}(t)|.
\label{est:psi-1}
\end{equation}

For the second term, we observe that also $s_{m,k}(t)$ lies in the subspace of $H$ generated by $e_{1}$, \ldots, $e_{k}$, and therefore from (\ref{th:snk}) with $m$ instead of $n$ we obtain that
\begin{equation}
|As_{m,k}(t)|\leq
\lambda_{k}|A^{1/2}s_{m,k}(t)|\leq
\lambda_{k}L_{1}^{1/2}.
\label{est:Asmk}
\end{equation}

Now we observe that
\begin{eqnarray*}
\lefteqn{\hspace{-1em}
|A^{1/2}u_{m}(t)|^{2}-|A^{1/2}u_{n}(t)|^{2}}
\\[0.5ex] 
& = &
|A^{1/2}s_{m,k}(t)|^{2}-|A^{1/2}s_{n,k}(t)|^{2}+|A^{1/2}r_{m,k}(t)|^{2}-|A^{1/2}r_{n,k}(t)|^{2}
\\[0.5ex] 
& = &
\langle A^{1/2}\rho_{n,m,k}(t),A^{1/2}(s_{m,k}(t)+s_{n,k}(t))\rangle
+|A^{1/2}r_{m,k}(t)|^{2}-|A^{1/2}r_{n,k}(t)|^{2},
\end{eqnarray*}
so that
\begin{eqnarray*}
\left||A^{1/2}u_{m}(t)|^{2}-|A^{1/2}u_{n}(t)|^{2}\right| & \leq &
|A^{1/2}\rho_{n,m,k}(t)|\left(|A^{1/2}s_{m,k}(t)|+|A^{1/2}s_{n,k}(t)|\right)
\\[0.5ex]
& &
\mbox{}+|A^{1/2}r_{m,k}(t)|^{2}+|A^{1/2}r_{n,k}(t)|^{2},
\end{eqnarray*}
and therefore from (\ref{th:snk}) and (\ref{th:rnk}) we obtain that
\begin{equation}
\left||A^{1/2}u_{m}(t)|^{2}-|A^{1/2}u_{n}(t)|^{2}\right|\leq
2L_{1}^{1/2}|A^{1/2}\rho_{n,m,k}(t)|+2R_{k}(u_{0},u_{1})L_{k}^{+}(t).
\label{est:Aum-Aun}
\end{equation}

From (\ref{hp:m-lip}), (\ref{est:Asmk}), and (\ref{est:Aum-Aun}) we deduce that
\begin{multline}
\left|m\left(|A^{1/2}u_{m}(t)|^{2}\right)-m\left(|A^{1/2}u_{n}(t)|^{2}\right)\right|
\cdot|As_{m,k}(t)|
\\[0.5ex]
\leq
2LL_{1}\lambda_{k}|A^{1/2}\rho_{n,m,k}(t)|+
2LL_{1}^{1/2}\lambda_{k}L_{k}^{+}(t)R_{k}(u_{0},u_{1}).
\label{est:psi-2}
\end{multline}

Plugging (\ref{est:psi-1}) and (\ref{est:psi-2}) into (\ref{psi-split}) we obtain that
\begin{eqnarray*}
|\psi_{n,k}(t)-\psi_{m,k}(t)| & \leq &
\lambda_{k}\left(1+\mu_{2}+2LL_{1}\right)|A^{1/2}\rho_{n,m,k}(t)|
\\[0.5ex]
& & 
\mbox{}+2LL_{1}^{1/2}\lambda_{k}L_{k}^{+}(t)R_{k}(u_{0},u_{1}),
\end{eqnarray*}
and therefore
\begin{eqnarray*}
G'(t) & \leq &
2|\psi_{n,k}(t)-\psi_{m,k}(t)|\cdot|\rho_{n,m,k}'(t)|
\\[0.5ex]
& \leq &
2\lambda_{k}\left(1+\mu_{2}+2LL_{1}\right)|A^{1/2}\rho_{n,m,k}(t)|\cdot|\rho_{n,m,k}'(t)|
\\
& &
\mbox{}+4LL_{1}^{1/2}\lambda_{k}L_{k}^{+}(t)R_{k}(u_{0},u_{1})|\rho_{n,m,k}'(t)|
\\[0.5ex]
& \leq &
\lambda_{k}\left(1+\mu_{2}+2LL_{1}\right)
\left(|\rho_{n,m,k}'(t)|^{2}+|A^{1/2}\rho_{n,m,k}(t)|^{2}\right)
\\
& &
\mbox{}+2LL_{1}\lambda_{k}|\rho_{n,m,k}'(t)|^{2}
+2L\lambda_{k}L_{k}^{+}(t)^{2}R_{k}(u_{0},u_{1})^{2}
\\[0.5ex]
& \leq &
\lambda_{k}\left(1+\mu_{2}+4LL_{1}\right)G(t)
+2L\lambda_{k}L_{k}^{+}(t)^{2}R_{k}(u_{0},u_{1})^{2}
\\[0.5ex]
& \leq &
\lambda_{k}\left(1+\mu_{2}+\frac{4LL_{1}}{\nu_{1}}\right)G(t)
+2L\lambda_{k}L_{k}^{+}(t)^{2}R_{k}(u_{0},u_{1})^{2}.
\end{eqnarray*}

Now let us set for simplicity
\begin{equation}
a_{k}:=
\left(1+\mu_{2}+\frac{4LL_{1}}{\nu_{1}}\right)\lambda_{k}=
\left(1+\mu_{2}+\frac{4LH_{0}}{\nu_{1}^{2}}\right)\lambda_{k}.
\label{defn:ak}
\end{equation}

Since $G(0)=0$, and $G(t)$ satisfies the differential inequality
\begin{equation}
G'(t)\leq a_{k}G(t)+2L\lambda_{k}L_{k}^{+}(t)^{2}R_{k}(u_{0},u_{1})^{2}
\qquad
\forall t\geq 0,
\nonumber
\end{equation}
we deduce that
\begin{equation}
G(t)\leq
2L\lambda_{k}R_{k}(u_{0},u_{1})^{2}
\exp(a_{k}t)\int_{0}^{t}L_{k}^{+}(s)^{2}\exp(-a_{k}s)\,ds
\qquad
\forall t\geq 0.
\nonumber
\end{equation}

On the other hand, from (\ref{defn:L+}) and (\ref{defn:ak}) we obtain that
\begin{eqnarray*}
\lefteqn{\hspace{-1em}
L_{k}^{+}(s)^{2}\exp(-a_{k}s)}
\\
& = & 
\left(\frac{\nu_{2}}{\nu_{1}}\right)^{2}
\exp\left(
\frac{2LH_{0}\lambda_{k}s}{\nu_{1}^{2}}
-(1+\mu_{2})\lambda_{k}s
-\frac{4LH_{0}\lambda_{k}s}{\nu_{1}^{2}}
\right)
\\
& \leq &
\left(\frac{\nu_{2}}{\nu_{1}}\right)^{2}
\exp(-(1+\mu_{2})\lambda_{k}s),
\end{eqnarray*}
and therefore
\begin{eqnarray*}
G(t) & \leq &
2L\lambda_{k}R_{k}(u_{0},u_{1})^{2}
\exp(a_{k}t)\left(\frac{\nu_{2}}{\nu_{1}}\right)^{2}
\int_{0}^{t}\exp(-(1+\mu_{2})\lambda_{k}s)\,ds
\\[0.5ex]
& \leq &
R_{k}(u_{0},u_{1})^{2}\cdot
2L\left(\frac{\nu_{2}}{\nu_{1}}\right)^{2}\exp(a_{k}t),
\end{eqnarray*}
which implies (\ref{th:s-s}).

\paragraph{\textmd{\textit{Conclusion}}}

Given initial data $(u_{0},u_{1})\in\SG(A)$, and given real numbers $T>0$ and $\ep>0$, we consider the constants $H_{0}$, $\mu_{1}$, $\mu_{2}$, $L$ that appear in (\ref{defn:E0}), (\ref{bound:m}) and (\ref{hp:m-lip}), and we choose a positive integer $k$ for which (\ref{hp:SG}) holds true, and such that
\begin{equation}
\max\left\{2L\left(\frac{\nu_{2}}{\nu_{1}}\right)^{2},
\frac{\nu_{2}}{\nu_{1}},
\left(1+\mu_{2}+\frac{4LH_{0}}{\nu_{1}^{2}}\right)T
\right\}\leq k,
\qquad\qquad
\ep\geq\frac{5}{k}.
\label{defn:k}
\end{equation}

For every pair of positive integers $(m,n)$, with $n>k$ and $m>k$, we write $u_{n}(t)$ and $u_{m}(t)$ as the sum of low-frequency and high-frequency components, and we obtain that
\begin{eqnarray}
\lefteqn{\hspace{-3em}
|u_{n}'(t)-u_{m}'(t)|^{2}+|A^{1/2}(u_{n}(t)-u_{m}(t))|^{2}}
\nonumber
\\[0.5ex]
& = &
|s_{n,k}'(t)-s_{m,k}'(t)|^{2}+|A^{1/2}(s_{n,k}(t)-s_{m,k}(t))|^{2}
\nonumber
\\[0.5ex]
& &
\mbox{}+|r_{n,k}'(t)-r_{m,k}'(t)|^{2}+|A^{1/2}(r_{n,k}(t)-r_{m,k}(t))|^{2}.
\label{un-um-split}
\end{eqnarray}

In the low-frequency regime we apply (\ref{th:s-s}). The first inequality in (\ref{defn:k}) implies that
$$L_{k}^{-}(t)\leq k\exp(k\lambda_{k})
\qquad
\forall t\in[0,T].$$
and therefore for every $t\in[0,T]$ it turns out that
\begin{eqnarray}
|s_{n,k}'(t)-s_{m,k}'(t)|^{2}+\left|A^{1/2}(s_{n,k}(t)-s_{m,k}(t))\right|^{2} & \leq & 
R_{k}(u_{0},u_{1})^{2}\cdot k\exp(k\lambda_{k})
\nonumber
\\[0.5ex]
& \leq & 
R_{k}(u_{0},u_{1})\cdot k\exp(k\lambda_{k}),
\label{est:low}
\end{eqnarray}
where the last inequality is true because (\ref{hp:SG}) implies in particular that $R_{k}(u_{0},u_{1})\leq 1$. 

In the high-frequency regime we estimate the norm of the difference with the sum of the norms, and from (\ref{th:rnk}) we obtain that
\begin{eqnarray*}
\lefteqn{\hspace{-2em}
|r_{n,k}'(t)-r_{m,k}'(t)|^{2}+\left|A^{1/2}(r_{n,k}(t)-r_{m,k}(t))\right|^{2}}
\\[0.5ex]
& \leq &
2\left(|r_{n,k}'(t)|^{2}+|A^{1/2}r_{n,k}(t)|^{2}\right)
+2\left(|r_{m,k}'(t)|^{2}+|A^{1/2}r_{m,k}(t)|^{2}\right)
\\[0.5ex]
& \leq &
4R_{k}(u_{0},u_{1})\cdot L_{k}^{+}(t).
\end{eqnarray*}

The first inequality in (\ref{defn:k}) implies that
$$L_{k}^{+}(t)\leq k\exp(k\lambda_{k})
\qquad
\forall t\in[0,T],$$
and therefore for every $t\in[0,T]$ it turns out that
\begin{equation}
|r_{n,k}'(t)-r_{m,k}'(t)|^{2}+\left|A^{1/2}(r_{n,k}(t)-r_{m,k}(t))\right|^{2}\leq 
R_{k}(u_{0},u_{1})\cdot 4k\exp(k\lambda_{k}).
\label{est:high}
\end{equation}

Plugging (\ref{est:low}) and (\ref{est:high}) into (\ref{un-um-split}) we conclude  that for every $t\in[0,T]$ it turns out that
\begin{equation}
|u_{n}'(t)-u_{m}'(t)|^{2}+|A^{1/2}(u_{n}(t)-u_{m}(t))|^{2}\leq
R_{k}(u_{0},u_{1})\cdot 5k\exp(k\lambda_{k})\leq
\frac{5}{k}\leq\ep,
\nonumber
\end{equation}
where in the second inequality we exploited again that (\ref{hp:SG}) is true for this value of $k$.

This completes the proof of the key claim (\ref{th:Cauchy}).
\qed


\setcounter{equation}{0}
\section{Final comments}\label{sec:final}

\paragraph{\textmd{\textit{Uniqueness}}}

Uniqueness of weak solutions is known only in the trivial case in which initial data have only a finite number of components. In this case indeed it is enough to consider the linear scalar problem (\ref{comp:eqn})--(\ref{comp:data}) satisfied by the components of $u(t)$, and deduce that components with zero initial data remain identically equal to zero for all positive times. What remains is a system of finitely many ordinary differential equations, for which uniqueness is a standard issue (if the nonlinearity $m(\sigma)$ is locally Lipschitz continuous, of course).

When initial data have infinitely many components different from zero, things are rather obscure. In particular, for regular initial data it is well-known that there exists a (local) strong solution, and this solution is unique \emph{in the class of strong solutions} (strong-strong uniqueness, see for example~\cite{1996-TAMS-AroPan} or the more recent~\cite{gg:K-uniqueness}), but we are not able to exclude the existence of a different weak solution, or even of infinitely many weak solutions, with the same initial condition. In other words, not only weak-weak uniqueness is an open problem, but also weak-strong uniqueness.

\paragraph{\textmd{\textit{Energy equality}}}

We do not know whether all weak solutions to (\ref{K:eqn}) satisfy the energy equality (\ref{eqn:energy-eq}). For sure it is true when the functions $t\mapsto|u'(t)|^{2}$ and $t\mapsto|A^{1/2}u(t)|^{2}$ are of class $C^{1}$, but this leads us once again to the space $D(A^{3/4})\times D(A^{1/4})$. Our solutions satisfy (\ref{eqn:energy-eq}) for a different reason, namely because they are uniform limits, with respect to the norm of the energy space, of strong solutions. On the other hand, nothing seems to prevent the existence of weak solutions, even with the same initial data, that do not satisfy the energy equality.

\paragraph{\textmd{\textit{Continuous dependence on initial data}}}

The limit of weak solutions (with respect to the norm of the energy space), if it exists, is again a weak solution. On the other hand, if a sequence of initial data $\{(u_{0,n},u_{1,n})\}$ converges in the energy space to some limit $(u_{0,\infty},u_{1,\infty})$, we do not know whether the corresponding sequence of weak solutions has a limit in the energy space, even up to subsequences.

The only case in which we have a positive answer is when $(u_{0,n},u_{1,n})\in\SG(A)$ and for every $n\geq 1$ the set of indices $k$ for which (\ref{hp:SG}) is true is the same (of course here we are considering ``our'' solutions, but we observed before that a priori different solutions might exist with the same initial data).


\subsubsection*{\centering Acknowledgments}

Both authors are members of the Italian {\selectlanguage{italian}%
``Gruppo Nazionale per l'Analisi Matematica, la Probabilit\`{a} e le loro Applicazioni'' (GNAMPA) of the ``Istituto Nazionale di Alta Matematica'' (INdAM)}. The first author was partially supported by PRIN 2020XB3EFL, ``Hamiltonian and Dispersive PDEs''. 

\selectlanguage{english}



\label{NumeroPagine}

\end{document}